# EXISTENCE, UNIQUENESS AND QUALITATIVE PROPERTIES OF POSITIVE SOLUTIONS OF QUASILINEAR ELLIPTIC EQUATIONS

PHUOC-TAI NGUYEN AND HOANG-HUNG VO

## Contents

ABSTRACT. We study the following quasilinear elliptic equation

$$-\Delta_p u + (\beta\Phi(x) - a(x))u^{p-1} + b(x)g(u) = 0 \quad \text{in } \mathbb{R}^N, \qquad (P_\beta)$$

where $p > 1$, $a, b \in L^\infty(\mathbb{R}^N)$, $\beta, b, g \geq 0$, $b \not\equiv 0$ and $\Phi \in L^\infty_{\text{loc}}(\mathbb{R}^N)$, $\inf_{\mathbb{R}^N} \Phi > -\infty$. We provide a sharp criterion in term of *generalized principal eigenvalues* for existence/nonexistence of positive solution of $(P_\beta)$ in suitable classes of functions. Uniqueness result for $(P_\beta)$ in those classes is also derived. Under additional conditions on $\Phi$, we further show that:

i) either for every $\beta \geq 0$ nonexistence phenomenon occurs,

ii) or there exists a threshold value $\beta^* > 0$ in the sense that for every $\beta \in [0, \beta^*)$ existence and uniqueness phenomenon occurs and for every $\beta \geq \beta^*$ nonexistence phenomenon occurs.

In the latter case, we study the limits, as $\beta \to 0$ and $\beta \to \beta^*$, of the sequence of positive solutions of $(P_\beta)$.

Our results are new even in the case $p = 2$.

## 1. Introduction and Main Results

1.1. **Introduction.** In this paper, we are interested in positive solutions of the quasilinear elliptic equation

$$(1.1) \qquad -\Delta_p u + (\beta\Phi(x) - a(x))u^{p-1} + b(x)g(u) = 0 \quad \text{in } \mathbb{R}^N,$$

where $\Delta_p u = \text{div}(|\nabla u|^{p-2}\nabla u)$, $p > 1$, $\Phi \in L^\infty_{\text{loc}}(\mathbb{R}^N)$, $\inf_{\mathbb{R}^N} \Phi > -\infty$, $a, b \in L^\infty(\mathbb{R}^N)$, $g \in C(\mathbb{R})$, $\beta, b, g \geq 0$, $b \not\equiv 0$.

By a (*weak*) *solution* of (1.1) we mean a nonnegative function $u \in W^{1,p}_{\text{loc}}(\mathbb{R}^N)$ such that $g(u) \in L^1_{\text{loc}}(\mathbb{R}^N)$ and $u$ satisfies (1.1) in the sense of distribution, namely

$$(1.2) \quad \int_{\mathbb{R}^N} (|\nabla u|^{p-2}\nabla u.\nabla\phi + (\beta\Phi - a)u^{p-1}\phi + bg(u)\phi)dx = 0 \qquad \forall\phi \in C^\infty_c(\mathbb{R}^N).$$

A nonnegative function $u$ is called a *supersolution* (resp. *subsolution*) of (1.1) if $u \in W^{1,p}_{\text{loc}}(\mathbb{R}^N)$, $g(u) \in L^1_{\text{loc}}(\mathbb{R}^N)$ and (1.2) holds with "=" replaced by "$\geq$" (resp. "$\leq$") and with nonnegative test function $\phi \in C^\infty_c(\mathbb{R}^N)$.

Equation (1.1) with $\beta = 0$ has a lot of applications in various aspects of mathematical biology, physics, especially in population dynamics. It has been intensively studied by many authors and numerous interesting results have been obtained in [2, 10, 6, 12, 14]. In particular, biological explanation of (1.1) was meticulously discussed in [6] for $p = 2$ and in [14] for $p > 1$. However, aside from [2, 3, 10], these papers concern only bounded or periodic domains and results involving unbounded domains (for instance $\mathbb{R}^N$) are much less. This paper is devoted to investigating a quasilinear version of the semilinear models in $\mathbb{R}^N$ proposed in [2, 3, 10]. Motivated by the above papers, we aim at establishing a criterion for the existence/nonexistence and uniqueness of positive solutions of (1.1). For this purpose, one of the main tasks is to study *the generalized principal eigenvalue* of the operator

$$(1.3) \qquad \mathcal{K}_V[\phi] := -\Delta_p\phi + V\phi^{p-1}, \quad \phi \geq 0, \qquad \text{in } \Omega,$$

where $\Omega \subset \mathbb{R}^N$ (possibly unbounded) and $V \in L^\infty_{\text{loc}}(\Omega)$, $\inf_\Omega V > -\infty$.

If $\Omega$ is a $C^{1,\theta}$ ($0 < \theta < 1$) bounded domain and $V \in L^\infty(\Omega)$, it is well-known that the variational problem

$$(1.4) \qquad \lambda_{1,V}(\Omega) := \inf_{\phi \in W^{1,p}_0(\Omega)\setminus\{0\}} \frac{\int_\Omega(|\nabla\phi|^p + V|\phi|^p)dx}{\int_\Omega |\phi|^p dx}$$

admits a unique (up to a scalar multiplication) minimizer $\varphi$ (see, e.g., [16, 21, Lemma 3]). Moreover, $\varphi$ is a $C^{1,\mu}$ ($0 < \mu < 1$) positive solution of the quasilinear eigenvalue problem

$$(1.5) \qquad \begin{cases} \mathcal{K}_V[\varphi] &= \lambda_{1,V}(\Omega)\varphi^{p-1} &\text{in } \Omega \\ \varphi &= 0 &\text{on } \partial\Omega. \end{cases}$$

$\lambda_{1,V}(\Omega)$ and $\varphi$ are called respectively the *principal eigenvalue* and *eigenfunction* of $\mathcal{K}_V$ in $\Omega$. Note that since $C_c^1(\Omega)$ is dense in $W_0^{1,p}(\Omega)$ with respect to $W^{1,p}$ norm, the infimum in (1.4) can be taken over $C_c^1(\Omega)$.

When $\Omega$ is an arbitrary (possibly unbounded) domain, following Berestycki et al. [7, 8, 9, 11, 28], we define

**Definition 1.1.** Put

$$\lambda(\mathcal{K}_V, \Omega) := \sup\{\lambda \in \mathbb{R} | \ \exists \psi \in W_{\mathrm{loc}}^{1,p}(\Omega), \psi > 0, \mathcal{K}_V[\psi] \geq \lambda \psi^{p-1} \ \ weakly \ \mathrm{in} \ \Omega\}.$$

$\lambda(\mathcal{K}_V, \Omega)$ is called the *generalized principal eigenvalue* of $\mathcal{K}_V$ in $\Omega$.

This type of eigenvalue was first introduced in a celebrated work of Berestycki-Nirenberg-Varadhan [8] for second order operators in bounded (not necessarily smooth) domains, and then was developed to second order operators in unbounded domains [7, 9, 11]. An important feature of the notion of generalized principal eigenvalue is that if $\Omega$ is a smooth and bounded domain, $\lambda(\mathcal{K}_V, \Omega)$ coincides with the principal eigenvalue $\lambda_{1,V}(\Omega)$, while if $\Omega$ is unbounded $\lambda(\mathcal{K}_V, \Omega)$ is well defined and can be expressed by a variational formula. For related definitions of generalized principal eigenvalues, the reader is referred to [7, 8, 9, 11] for linear operators, [4, 28] for fully nonlinear operators and [13] for singular fully nonlinear operators. To our knowledge, no investigation of generalized principal eigenvalue for quasilinear operators has been previously obtained. A related notion, which is called *the best constant in the Hardy-type inequality*, was used by Pinchover to study optimal Hardy-type inequalities (see [17]). Our approach is different from the quoted works and the main contributions are:

- the investigation of the generalized principal eigenvalue of $\mathcal{K}_V$,
- the existence/nonexistence and uniqueness of positive solution of (1.1),
- the study of threshold value for existence of positive solution of (1.1).

We say that $\{\Omega_n\}$ is an *exhaustion* of $\Omega$ if $\Omega_n$ is a $C^{1,\mu}$ ($\mu$ may depend on $n$) bounded domain, $\Omega_n \Subset \Omega_{n+1}$ and $\cup_n \Omega_n = \Omega$.

We first present results concerning qualitative properties of the generalized principal eigenvalues.

**Theorem 1.2.** (1) Assume $\Omega$ is a $C^{1,\theta}$ ($\theta \in (0,1)$) bounded domain in $\mathbb{R}^N$ and $V \in L^\infty(\Omega)$. Then

$$\lambda(\mathcal{K}_V, \Omega) = \lambda_{1,V}(\Omega).$$

(2) Assume $\Omega$ is a domain in $\mathbb{R}^N$ (possibly unbounded) with $\{\Omega_n\}$ is an exhaustion of $\Omega$. Let $V \in L_{\mathrm{loc}}^\infty(\Omega)$ and $\inf_\Omega V > -\infty$. Then the following properties hold.

(i) $\inf_\Omega V \leq \lambda(\mathcal{K}_V, \Omega_{n+1}) < \lambda(\mathcal{K}_V, \Omega_n)$ for every $n \in \mathbb{N}$.

(ii) $\lambda(\mathcal{K}_V, \Omega) = \lim_{n \to \infty} \lambda(\mathcal{K}_V, \Omega_n)$ and there exists a positive weak solution $\varphi \in C_{\mathrm{loc}}^1(\Omega)$ of

$$\mathcal{K}_V[\varphi] = \lambda(\mathcal{K}_V, \Omega)\varphi^{p-1} \qquad \mathrm{in} \ \Omega.$$

(iii)

$$\lambda(\mathcal{K}_V, \Omega) = \inf_{\phi \in C_c^1(\Omega) \setminus \{0\}} \frac{\int_\Omega (|\nabla \phi|^p + V|\phi|^p)dx}{\int_\Omega |\phi|^p dx}. \tag{1.6}$$

The notion of generalized principal eigenvalue plays a crucial role in proving the existence of positive solutions of (1.1). Before stating the next results, let us describe the hypotheses that we need in the paper.

(A1) There exist $\alpha \in [0, p]$ and $m > 0$ such that $\limsup_{|x| \to \infty} |x|^\alpha a(x) < -m$.

(A2) The function $\frac{g(s)}{s^{p-1}}$ is increasing on $(0, \infty)$ and $\lim_{s \to 0} \frac{g(s)}{s^{p-1}} = 0$.

(A3) There exists $s_0 > 0$ such that $-a(x)s_0^{p-1} + b(x)g(s_0) \geq 0$ a.e. in $\mathbb{R}^N$.

(A4) The set $\{x \in \mathbb{R}^N : \Phi(x) < 0\}$ is bounded.

All above assumptions have biological interpretations. Hypothesis (A1) refers to the environments being unfavorable, unfavorably neutral or nearly neutral near infinity, according to the cases $\alpha = 0$, $\alpha \in (0, p)$ or $\alpha = p$. This kind of assumption is recently used to describe the effect of global warming (see [2],[5],[10],[32]). Hypothesis (A2) means that the intrinsic growth rate decreases when the population density is increasing. This is due to the intraspecific competition for resources. Hypothesis (A3) expresses a saturation effect: when the population density is large, the death rate is higher than the birth rate and the population decreases. Lastly, (A4) allows a refuge zone so that the population may persist. This zone is assumed to be bounded and surrounded by a hostile environment.

As it will be shown below, the decay rate of the potential $a$ stated in (A1) has a significant effect on solutions of (1.1). To be precise, solutions of (1.1) decay exponentially when $a$ is a slow decay potential (i.e. $\alpha \in [0, p)$), while they decay polynomially when $a$ is a Hardy potential (i.e. $\alpha = p$). As an example, one can take $a(x)$ to behave exactly like $-m/2|x|^{-\alpha}$ as $|x| \to \infty$.

In the sequel, when $\Omega = \mathbb{R}^N$ and $V = \beta\Phi - a$ with $\beta \geq 0$, for simplicity and to emphasize the dependence of $\mathcal{K}$ on $\beta$, we use the notation $\mathcal{K}_\beta$ in stead of $\mathcal{K}_V$, i.e.

$$\mathcal{K}_\beta[\phi] := -\Delta_p \phi + (\beta\Phi - a)\phi^{p-1}, \quad \phi \geq 0. \tag{1.7}$$

Denote by $\lambda_\beta$ the generalized principal eigenvalue of $\mathcal{K}_\beta$, i.e. $\lambda_\beta = \lambda(\mathcal{K}_\beta, \mathbb{R}^N)$. In this case, we are able to obtain the existence, nonexistence and uniqueness result for (1.1) in some classes of functions. To this aim, let us define these classes.

For $p \geq 2$ and $\alpha \in [0, p)$, set

$$\mathcal{S}_p = \left\{ u \mid u > 0 \text{ in } \mathbb{R}^N \text{ and } \limsup_{|x| \to \infty} e^{-\overline{\omega}|x|^{1-\frac{\alpha}{p}}} u(x) < \infty \right\}$$

where

$$(1.8) \qquad \overline{\omega} = \left(\frac{2m}{N(p-1)}\right)^{1/p} \frac{p}{p-\alpha}$$

with $m$ being in (A1).

For $p \geq 2$ and $\alpha = p$, set

$$\mathcal{S}_p = \{u \mid u > 0 \text{ in } \mathbb{R}^N \text{ and } \limsup_{|x| \to \infty} |x|^{-\gamma_0} u(x) < \infty\}$$

where $\gamma_0$ is the unique positive solution of algebraic equation

$$(1.9) \qquad ((p-2)N+1)\gamma^p + (3p-5)N\gamma^{p-1} - m = 0.$$

One of the main results of the paper is the following.

**Theorem 1.3.** Assume $p \geq 2$, $\alpha \in [0, p]$, $\beta \geq 0$ and hypotheses (A1)-(A4) are satisfied.

(1) If $\lambda_\beta < 0$ then there exists a unique positive solution $u_\beta$ of (1.1) in $\mathcal{S}_p$. Moreover,

  (i) If $\alpha \in [0, p)$ then

$$(1.10) \qquad \lim_{|x| \to \infty} e^{\overline{\omega}|x|^{1-\frac{\alpha}{p}}} u_\beta(x) = 0.$$

  (ii) If $\alpha = p$ then

$$(1.11) \qquad \lim_{|x| \to \infty} |x|^{\gamma_0} u_\beta(x) = 0.$$

(2) If $\lambda_\beta \geq 0$ then there exists no positive solution of (1.1) in $\mathcal{S}_p$.

When $p \geq 2$, due to a-priori estimates (see Section 4.2), solutions of (1.1) belonging to $\mathcal{S}_p$ decay. Note that $\mathcal{S}_p$ is *larger than the class of bounded functions*. This, together with the comparison principle (Theorem 3.1), implies the uniqueness in $\mathcal{S}_p$. When $1 < p < 2$, because of the lack of such a-priori estimates, uniqueness might not hold in $\mathcal{S}_p$, or even in the class of bounded positive solutions. Therefore, we only deal with decaying solutions of (1.1). Thus, when $1 < p < 2$ the class of functions under consideration is

$$\mathcal{S}_p = \{u \mid u > 0 \text{ in } \mathbb{R}^N \text{ and } \lim_{|x| \to \infty} u(x) = 0\}.$$

**Theorem 1.4.** Assume $1 < p < 2$, $\alpha \in [0, p]$, $\beta \geq 0$ and hypotheses (A1)-(A4) are satisfied.

(1) If $\lambda_\beta < 0$ then there exists a unique solution $u_\beta$ of (1.1) in $\mathcal{S}_p$. Moreover,

  (i) If $\alpha \in [0, p)$ then $u_\beta$ decays exponentially.

  (ii) If $\alpha = p$ then $u_\beta$ decays polynomially.

(2) If $\lambda_\beta \geq 0$ then there exists no solution of (1.1) in $\mathcal{S}_p$.

**Remark 1.** Such type of results might be also called Liouville type result for (1.1). For other types of Liouville results related to semilinear and quasilinear equations, interested readers are referred to [7, 18, 27].

In order to state the next result, we introduce the following weighted spaces. Assume that $\Phi \geq 0$, $\Phi \in L^\infty_{\mathrm{loc}}(\mathbb{R}^N)$. For $p \geq 1$, let $D^{1,p}(\mathbb{R}^N)$ be the completion of $C^\infty_c(\mathbb{R}^N)$ with respect to the norm

$$\|u\|_{D^{1,p}(\mathbb{R}^N)} = \left( \int_{\mathbb{R}^N} |\nabla u|^p dx \right)^{1/p}.$$

Denote $L^p_\Phi(\mathbb{R}^N) = \{u : \mathbb{R}^N \to \mathbb{R} \text{ measurable } : \int_{\mathbb{R}^N} \Phi |u|^p dx < \infty\}$ with the norm

$$\|u\|_{L^p_\Phi(\mathbb{R}^N)} = \left( \int_{\mathbb{R}^N} \Phi |u|^p dx \right)^{1/p}.$$

Define $W^{1,p}_\Phi(\mathbb{R}^N) = D^{1,p}(\mathbb{R}^N) \cap L^p_\Phi(\mathbb{R}^N)$ then $W^{1,p}_\Phi(\mathbb{R}^N)$ is a Banach space with respect to the norm

$$\|u\|_{W^{1,p}_\Phi(\mathbb{R}^N)} = \left( \int_{\mathbb{R}^N} (|\nabla u|^p + \Phi |u|^p) dx \right)^{1/p}.$$

We assume that

(A5)     $\Phi(0) \geq 0$, $\Phi > 0$ in $\mathbb{R}^N \setminus \{0\}$, $\Phi \in L^\infty_{\mathrm{loc}}(\mathbb{R}^N)$, $\lim_{|x| \to \infty} \Phi(x) = \infty$ and
          the embedding $W^{1,p}_\Phi(\mathbb{R}^N) \hookrightarrow L^p(\mathbb{R}^N)$ is compact.

**Remark 2.** When $p = 2$ the embedding $W^{1,p}_\Phi(\mathbb{R}^N) \hookrightarrow L^p(\mathbb{R}^N)$ is compact only if $\lim_{|x| \to \infty} \Phi(x) = \infty$. When $p > 1$, the compact embedding holds, for instance, if $1 < p < N$, $\Phi^{-\frac{\alpha_1 \ell}{(1-\alpha_1)(\ell-p)}} \in L^1(\mathbb{R}^N)$ and $\Phi^{-\frac{\alpha_2 N}{(1-\alpha_2)p}} \in L^1(\mathbb{R}^N)$ with $\alpha_1, \alpha_2 \in (0,1)$, $\ell \in (p, \frac{Np}{N-p})$ (See [30, 33] for more detail).

Our last result concerns the existence of a threshold value for existence of positive solutions of (1.1).

**Theorem 1.5.** Assume $p > 1$, $\alpha \in [0, p]$ and hypotheses (A1), (A2), (A3) and (A5) are satisfied. Let $\lambda_0$ be the generalized principal eigenvalue of $\mathcal{K}_0$ in $\mathbb{R}^N$.

(1) If $\lambda_0 \geq 0$ then for every $\beta > 0$ there is no positive solution of (1.1) in $\mathcal{S}_p$.

(2) If $\lambda_0 < 0$ then there exists a threshold value $\beta^*$ in the following sense:

  (i) For every $\beta \in [0, \beta^*)$ there exists a unique solution $u_\beta$ of (1.1) in $\mathcal{S}_p$.

  (ii) For every $\beta \in [\beta^*, \infty)$ there is no solution of (1.1) in $\mathcal{S}_p$.

Moreover $u_\beta \to u_0$ as $\beta \to 0$ and $u_\beta \to 0$ as $\beta \to \beta^*$ in $L^q(\mathbb{R}^N)$ for all $q \in (0, \infty]$.

1.2. **Organization of the paper and strategy of the proofs.** The rest of the paper is organized as follows.

Section 2 is devoted to the study of the generalized principal eigenvalue and its properties in bounded and unbounded domains.

In Section 3, we establish a comparison principle for (1.1) in $\mathbb{R}^N$. Although the technique is inspired from [16, 20], we have to overcome extra difficulties stemming from the lack of compactness in domain. It is worth emphasizing that since $\Delta_p(u + v) \neq \Delta_p u + \Delta_p v$, the sliding argument based on strong maximum principle and the variational argument used for semilinear equations in [2, 10] fails to apply to this framework. Our comparison principle allows to compare a sub and a supersolution of (1.1) without extra assumption on the gradient of supersolution as in [10, Theorem 2.3].

In Section 4, we prove Theorem 1.3 and Theorem 1.4. Existence part is demonstrated by sub and super solution method. When $1 < p < 2$, in order to obtain nonexistence and uniqueness, we use the comparison principle. When $p \geq 2$, we solve equation (1.1) in a larger class of solutions, which may be unbounded in $\mathbb{R}^N$. Due to a delicate construction of supersolutions in exterior domains, we show that all positive solutions of (1.1) in such a class decay. The rate of decay in two cases $\alpha \in [0, p)$ and $\alpha = p$ are different. More precisely, they decay exponentially in the first case while they decay polynomially in the second case. Our construction relies essentially on the weak comparison principle (see [15]) and a-priori growth condition. Notice that, throughout the construction, boundedness assumption of solution is relaxed.

Finally, in Section 5, thanks to the compact embedding provided in (A5), we prove Theorem 1.5.

## 2. Generalized principal eigenvalue

This section is devoted to the investigation of the generalized principal eigenvalue $\lambda(\mathcal{K}_V, \Omega)$ of $\mathcal{K}_V$. Denote

$$\Lambda_V(\Omega) := \{\lambda \in \mathbb{R} |\ \exists \psi \in W^{1,p}_{\text{loc}}(\Omega), \psi > 0, \mathcal{K}_V[\psi] \geq \lambda \psi^{p-1}\ \text{ weakly in } \Omega\}.$$

Note that $\Lambda_V(\Omega) \neq \emptyset$ since $\inf_{\mathbb{R}^N} V > -\infty$.

**Lemma 2.1.** Assume $\Omega \subset \mathbb{R}^N$ is a $C^{1,\theta}$ bounded domain for some $\theta \in (0, 1)$ and $V \in L^\infty(\Omega)$. Let $\lambda_{1,V}$ be the principal eigenvalue of $\mathcal{K}_V$ given in (1.4). Then $\lambda \leq \lambda_{1,V}(\Omega)$ for every $\lambda \in \Lambda_V(\Omega)$.

*Proof.* Let $\lambda \in \Lambda_V(\Omega)$ then there exists $\psi > 0$, $\psi \in W^{1,p}_{\text{loc}}(\Omega)$ such that

$$(2.1) \qquad \int_\Omega |\nabla \psi|^{p-2} \nabla \psi \cdot \nabla \phi dx + \int_\Omega V \psi^{p-1} \phi dx \geq \lambda \int_\Omega \psi^{p-1} \phi dx \quad \forall \phi \in C_c^\infty(\Omega).$$

Let $\varphi_V^\Omega$ be the eigenfunction associated with $\lambda_{1,V}(\Omega)$, normalized by $\left\|\varphi_V^\Omega\right\|_{L^p(\Omega)} = 1$. It is classical that $\varphi_V^\Omega \in W_0^{1,p}(\Omega)$ and

$$\text{(2.2)} \qquad \int_\Omega |\nabla\varphi_V^\Omega|^p dx + \int_\Omega V(\varphi_V^\Omega)^p dx = \lambda_{1,V}(\Omega).$$

Since $\varphi_V^\Omega \in W_0^{1,p}(\Omega)$, one can find a sequence $\phi_n \in C_c^\infty(\Omega)$ such that $\phi_n \to \varphi_V^\Omega$ in $W^{1,p}(\Omega)$. Then $\frac{\phi_n^p}{\psi^{p-1}}$ can be used as a test function in (2.1) to obtain

$$\text{(2.3)} \qquad \int_\Omega |\nabla\psi|^{p-2}\nabla\psi \cdot \nabla\left(\frac{\phi_n^p}{\psi^{p-1}}\right) dx + \int_\Omega V\phi_n^p dx \geq \lambda \int_\Omega \phi_n^p dx.$$

Subtracting (2.3) from (2.2) yields

$$\text{(2.4)} \qquad \begin{aligned} &\lambda_{1,V}(\Omega) - \lambda\int_\Omega \phi_n^p dx \\ &\geq \int_\Omega (|\nabla\varphi_V^\Omega|^p - |\nabla\phi_n|^p) dx + \int_\Omega \left[|\nabla\phi_n|^p - |\nabla\psi|^{p-2}\nabla\psi \cdot \nabla\left(\frac{\phi_n^p}{\psi^{p-1}}\right)\right] dx \\ &\quad + \int_\Omega V[(\varphi_V^\Omega)^p - \phi_n^p] dx \\ &:= I_{n,1} + I_{n,2} + I_{n,3}. \end{aligned}$$

Applying Picone's identity [1, Theorem 1.1]

$$\begin{aligned} |\nabla\xi_1|^p &- |\nabla\xi_2|^{p-2}\nabla\xi_2\nabla\left(\frac{\xi_1^p}{\xi_2^{p-1}}\right) \\ &= |\nabla\xi_1|^p + (p-1)\frac{\xi_1^p}{\xi_2^p}|\nabla\xi_2|^p - p\frac{\xi_1^{p-1}}{\xi_2^{p-1}}\nabla\xi_1|\nabla\xi_2|^{p-2}\nabla\xi_2 \geq 0, \end{aligned}$$

with $\xi_1 = \phi_n$ and $\xi_2 = \psi$, we deduce $I_{n,2} \geq 0$. In fact, Picone's identity in [1] is proved for $C^1$ functions $\xi_1$, $\xi_2$, however, one can show that Picone's identity remains true if $\xi_1, \xi_2 \in W^{1,p}(\Omega)$. Since $V \in L^\infty(\Omega)$ and $\phi_n \to \varphi_V^\Omega$ in $W^{1,p}(\Omega)$, $\lim_{n\to\infty} I_{n,1} = \lim_{n\to\infty} I_{n,3} = 0$. Therefore, by letting $n \to \infty$ in (2.4), we obtain $\lambda_{1,V}(\Omega) - \lambda \geq 0$. □

**Remark 3.** This lemma may be obtained in a different way due to the aid of a general version of Allegretto-Piepenbrink theorem (for more details, see [26, Theorem 2.3]).

**Proof of Theorem 1.2.** Statement (1) follows from Lemma 2.1 and the fact that $\lambda_{1,V}(\Omega) \leq \lambda(\mathcal{K}_V, \Omega)$.

Statement (2.i) follows from Statement (1) and the strict monotonicity property of the principal eigenvalue.

We next prove statement (2.ii). From (2.i), there exists $\bar\lambda_V := \lim_{n\to\infty} \lambda(\mathcal{K}_V, \Omega_n)$ and $\bar\lambda_V \geq \lambda(\mathcal{K}_V, \Omega)$. We will show that $\bar\lambda_V = \lambda(\mathcal{K}_V, \Omega)$. Fix $x_0 \in \Omega$ such that $x_0 \in \Omega_n$ for every $n \in \mathbb{N}$. Due to statement (1), for every $n \in \mathbb{N}$, $\lambda(\mathcal{K}_V, \Omega_n)$ can be variationally characterized by

$$(2.5) \qquad \lambda(\mathcal{K}_V, \Omega_n) = \inf_{\phi \in C_c^1(\Omega_n) \setminus \{0\}} \frac{\int_{\Omega_n} (|\nabla \phi|^p + V|\phi|^p) dx}{\int_{\Omega_n} |\phi|^p dx}.$$

Let $\varphi_V^{\Omega_n}$ be the eigenfunction associated with $\lambda(\mathcal{K}_V, \Omega_n)$, normalized by $\varphi_V^{\Omega_n}(x_0) = 1$. Then $\varphi_V^{\Omega_n} \in W_0^{1,p}(\Omega_n)$ and $\varphi_V^{\Omega_n}$ is a weak solution of

$$\mathcal{K}_V[\varphi] = \lambda(\mathcal{K}_V, \Omega_n) \varphi^{p-1}$$

in $\Omega_n$. Take an arbitrary subset $G \Subset \Omega$ then there exists $n_G > 0$ such that $G \Subset \Omega_n$ for every $n \geq n_G$. By Harnack's inequality (see, e.g., [29, Theorem 5], [31, Theorem 1]), one can find a constant $C_G > 0$ independent of $n$ such that

$$\sup_G \varphi_V^{\Omega_n} < C_G \inf_G \varphi_V^{\Omega_n} \leq C_G \quad \forall n \geq n_G.$$

In light of local regularity for elliptic equations and standard argument, up to a subsequence, $\{\varphi_V^{\Omega_n}\}$ converges in $C^1_{\mathrm{loc}}(\Omega)$ to a function $\overline{\varphi}$ which is a weak solution of

$$\mathcal{K}_V[\varphi] = \bar{\lambda}_V \varphi^{p-1}$$

in $\Omega$. Moreover, $\overline{\varphi}(x_0) = 1$ and therefore $\overline{\varphi} > 0$ in $\Omega$ by Harnack's inequality. Therefore $\bar{\lambda}_V \in \Lambda_V(\Omega)$ and consequently $\bar{\lambda}_V \leq \lambda(\mathcal{K}_V, \Omega)$. Thus $\bar{\lambda}_V = \lambda(\mathcal{K}_V, \Omega)$.

Finally, we prove (2.iii). Denote by $\tilde{\lambda}_V$ the right-hand side of (1.6) then $\tilde{\lambda}_V \geq \inf_\Omega V > -\infty$. Thanks to (2.5), $\lambda(\mathcal{K}_V, \Omega_n) \geq \tilde{\lambda}_V$ for every $n \in \mathbb{N}$ and consequently $\lambda(\mathcal{K}_V, \Omega) \geq \tilde{\lambda}_V$. We next prove the inverse inequality. By the definition of $\tilde{\lambda}_V$, for any $\delta > 0$, there exists $\phi_\delta \in C_c^1(\Omega)$, $\phi_\delta \neq 0$ such that

$$(2.6) \qquad \tilde{\lambda}_V + \delta \geq \frac{\int_\Omega (|\nabla \phi_\delta|^p + V|\phi_\delta|^p) dx}{\int_\Omega |\phi_\delta|^p dx}.$$

Since $\phi_\delta \in C_c^1(\Omega)$, there exists $N_\delta > 0$ such that $\phi_\delta \in C_c^1(\Omega_n)$ for every $n \geq N_\delta$. Due to (2.5), the right-hand side of (2.6) is greater than $\lambda(\mathcal{K}_V, \Omega_n)$ for every $n \geq N_\delta$. Consequently

$$\tilde{\lambda}_V + \delta \geq \lambda(\mathcal{K}_V, \Omega_n) \quad \forall n \geq N_\delta.$$

By letting $n \to \infty$ and $\delta \to 0$ successively, we obtain $\tilde{\lambda}_V \geq \lambda(\mathcal{K}_V, \Omega_n)$ and therefore $\tilde{\lambda}_V = \lambda(\mathcal{K}_V, \Omega)$. $\qquad\square$

The next result concerns the case $\Omega = \mathbb{R}^N$ and $V = \beta\Phi - a$. In this case, to simplify notations and to to emphasize the dependence of $\mathcal{K}$ on $\beta$, we denote $\mathcal{K}_\beta$ instead of $\mathcal{K}_V$, i.e. $\mathcal{K}_\beta[\phi] := -\Delta_p \phi + (\beta\Phi - a)\phi^{p-1}$.

**Proposition 2.2.** Assume $\Omega = \mathbb{R}^N$ and $V = \beta\Phi - a$ with $\beta \geq 0$, $\Phi \in L_{\mathrm{loc}}^\infty(\mathbb{R}^N)$, $\Phi \geq 0$ and $a \in L^\infty(\mathbb{R}^N)$. Let $\lambda_\beta$ be the generalized principal eigenvalue of $\mathcal{K}_\beta$ defined in (1.7). Then $\lim_{\beta \to 0} \lambda_\beta = \lambda_0$.

*Proof.* Notice that $\lambda(\mathcal{K}_{\beta'}, \Omega) \geq \lambda(\mathcal{K}_\beta, \Omega) \geq \lambda(\mathcal{K}_0, \Omega)$ for every $0 < \beta < \beta'$ since $\Phi \geq 0$. Then there exists $\underline{\lambda} := \lim_{\beta \to 0} \lambda(\mathcal{K}_\beta, \Omega)$ and $\underline{\lambda} \geq \lambda(\mathcal{K}_0, \Omega)$. By variational characterization of $\lambda(\mathcal{K}_0, \Omega)$, for any $\delta > 0$, there exists $\phi_\delta \in C_c^1(\Omega)$, $\|\phi_\delta\|_{L^p(\Omega)} = 1$, such that

$$\lambda(\mathcal{K}_0, \Omega) + \delta \geq \int_\Omega (|\nabla \phi_\delta|^p - a|\phi_\delta|^p) dx.$$

Since $\Phi \in L_{\text{loc}}^\infty(\Omega), \phi_\delta \in C_c^1(\Omega)$, we get $\Phi|\phi_\delta|^p \in L^1(\Omega)$. Choosing $\beta$ such that $\beta \int_\Omega \Phi|\phi_\delta|^p < \delta$, we get

$$\lambda(\mathcal{K}_0, \Omega) + 2\delta \geq \int_\Omega [|\nabla \phi_\delta|^p + (\beta\Phi - a)|\phi_\delta|^p] dx \geq \lambda(\mathcal{K}_\beta, \Omega) \geq \underline{\lambda}.$$

Letting $\delta \to 0$ yields $\lambda(\mathcal{K}_0, \Omega) \geq \underline{\lambda}$, which concludes the proof. $\square$

## 3. Comparison principle

In this section, we prove the comparison principle for (1.1) which serve to obtain the nonexistence and uniqueness result in the next section. It is noteworthy that the comparison principle is never obvious for quasilinear operators. In the sequel, $B_R$ denotes the ball of radius $R > 0$ and center $0$.

**Theorem 3.1.** Let $p > 1$, $\beta \geq 0$, $\Phi \in L_{\text{loc}}^\infty(\mathbb{R}^N)$, $a, b \in L^\infty(\mathbb{R}^N)$, $b \geq 0$, $b \not\equiv 0$ and (A2) hold. Assume $\mathcal{N} := \{x \in \mathbb{R}^N : a(x) - \beta\Phi(x) \geq 0\}$ is bounded. Let $u_1, u_2 \in C_{\text{loc}}^1(\mathbb{R}^N)$ be respectively positive supersolution and subsolution of (1.1) such that

$$(3.1) \qquad\qquad \lim_{|x| \to \infty} u_2(x) = 0.$$

Then $u_1 \geq u_2$ in $\mathbb{R}^N$.

**Remark 4.** In Theorem 3.1 we do not require $a$ to satisfy (A1) and $\Phi$ to be non-negative. It is clear that $\mathcal{N}$ is bounded if (A1) is fulfilled and $\Phi \geq 0$.

*Proof.* Since $b \not\equiv 0$ in $\mathbb{R}^N$, there exists $R_b > 0$ such that $b \not\equiv 0$ in $B_{R_b}$. For any $\epsilon \geq 0$, $R \in (0, \infty]$, we set

$$D_R(\epsilon) = \{x \in B_R : u_2(x) + \epsilon > u_1(x) + 2\epsilon\}.$$

We will show that $D_\infty(0) = \emptyset$. It suffices to prove that $D_{R_0}(0) = \emptyset$ for all $R_0 > R_b$.

Fix $R_0 > R_b$ and let $\epsilon \in (0, 1)$. By (3.1), there exists $R = R(\epsilon) > R_0$ such that $u_2(x) < \epsilon/2$ for every $x \in B_R^c$. Denote $\epsilon_1 = 2\epsilon$, $\epsilon_2 = \epsilon$ and

$$v_i = \frac{[(u_2 + \epsilon_2)^p - (u_1 + \epsilon_1)^p]^+}{(u_i + \epsilon_i)^{p-1}} \quad i = 1, 2,$$

where $h^+ = \max\{h, 0\}$. Then $v_i \in W_0^{1,p}(K_R)$ with some $K_R \Subset B_R$ and $v_i = 0$ outside $K_R$. Therefore $v_i$ can be approximated by a sequence of functions in $C_0^\infty(K_R)$. We can use $v_i$, $i = 1, 2$, as test functions to get

$$\tag{3.2} \int_{B_R} (|\nabla u_1|^{p-2}\nabla u_1.\nabla v_1 + (\beta\Phi - a)u_1^{p-1}v_1 + bg(u_1)v_1)dx \geq 0,$$

$$\tag{3.3} \int_{B_R} (|\nabla u_2|^{p-2}\nabla u_2.\nabla v_2 + (\beta\Phi - a)u_2^{p-1}v_2 + bg(u_2)v_2)dx \leq 0.$$

Subtracting (3.3) from (3.2) yields

$$\tag{3.4} \begin{aligned} \int_{B_R} (|\nabla u_1|^{p-2}\nabla u_1.\nabla v_1 &- |\nabla u_2|^{p-2}\nabla u_2.\nabla v_2)dx \\ &\geq -\int_{B_R} (\beta\Phi - a)(u_1^{p-1}v_1 - u_2^{p-1}v_2)dx - \int_{B_R} b(g(u_1)v_1 - g(u_2)v_2)dx. \end{aligned}$$

Set

$$\tag{3.5} w_i = u_i + \epsilon_i, \quad W_i = \nabla(\ln(u_i + \epsilon_i)) = \frac{\nabla u_i}{u_i + \epsilon_i} \quad i = 1, 2$$

and

$$I := |\nabla u_2|^{p-2}\nabla u_2 \cdot \nabla v_2 - |\nabla u_1|^{p-2}\nabla u_1 \cdot \nabla v_1.$$

By a computation, we obtain

$$\tag{3.6} \begin{aligned} I = \ & w_2^p(|W_2|^p - |W_1|^p - p|W_1|^{p-2}W_1 \cdot (W_2 - W_1)) \\ & + w_1^p(|W_1|^p - |W_2|^p - p|W_2|^{p-2}W_2 \cdot (W_1 - W_2)). \end{aligned}$$

By Clarkson's inequality [25], for all vectors $X, Y \in \mathbb{R}^N$, we have

$$\tag{3.7} |X|^p - |Y|^p - p|Y|^{p-2}Y(X - Y) \geq c_p \frac{|X - Y|^{p+(2-p)^+}}{(|X| + |Y|)^{(2-p)^+}}$$

where

$$c_p = \begin{cases} \dfrac{1}{2^p - 1} & \text{if } 1 < p < 2 \\ \dfrac{3p(p-1)}{16} & \text{if } p \geq 2. \end{cases}$$

It follows from (3.6) and (3.7) that

$$\tag{3.8} I \geq c_p(w_1^p + w_2^p)\frac{|W_1 - W_2|^{p+(2-p)^+}}{(|W_1| + |W_2|)^{(2-p)^+}} \quad \forall p > 1.$$

Combining (3.4) and (3.8) and using (3.5), we obtain

$$
\begin{aligned}
-c_p \int_{D_R(\epsilon)} (w_1^p + w_2^p) &\frac{|W_1 - W_2|^{p+(2-p)^+}}{(|W_1| + |W_2|)^{(2-p)^+}} dx \\
&\geq \int_{D_R(\epsilon)} (\beta\Phi - a)\left[\left(\frac{u_2}{w_2}\right)^{p-1} - \left(\frac{u_1}{w_1}\right)^{p-1}\right](w_2^p - w_1^p)dx \\
&\quad + \int_{D_R(\epsilon)} b\left[\frac{g(u_2)}{w_2^{p-1}} - \frac{g(u_1)}{w_1^{p-1}}\right](w_2^p - w_1^p)dx \\
&\geq \int_{D_R(\epsilon)\cap\mathcal{N}} (\beta\Phi - a)\left[\left(\frac{u_2}{w_2}\right)^{p-1} - \left(\frac{u_1}{w_1}\right)^{p-1}\right](w_2^p - w_1^p)dx \\
&\quad + \int_{D_R(\epsilon)} b\left[\frac{g(u_2)}{w_2^{p-1}} - \frac{g(u_1)}{w_1^{p-1}}\right](w_2^p - w_1^p)dx
\end{aligned}
\tag{3.9}
$$

where $\mathcal{N} = \{x \in \mathbb{R}^N : a(x) - \beta\Phi(x) \geq 0\}$. The last inequality in (3.9) is derived from the fact that

$$
\left[\left(\frac{u_2}{w_2}\right)^{p-1} - \left(\frac{u_1}{w_1}\right)^{p-1}\right](w_2^p - w_1^p) \geq 0 \quad \text{in } D_R(\epsilon).
$$

Since $\mathcal{N}$ is bounded, we derive from (3.9) that

$$
\begin{aligned}
c_p \int_{D_R(\epsilon)} (w_1^p + w_2^p)&\frac{|W_2 - W_1|^{p+(2-p)^+}}{(|W_1| + |W_2|)^{(2-p)^+}} dx + \int_{D_R(\epsilon)} b\left[\frac{g(u_2)}{w_2^{p-1}} - \frac{g(u_1)}{w_1^{p-1}}\right](w_2^p - w_1^p)dx \\
&\leq \sup_{\mathcal{N}} |a - \beta\Phi| \int_{D_R(\epsilon)\cap\mathcal{N}} \left[\left(\frac{u_2}{w_2}\right)^{p-1} - \left(\frac{u_1}{w_1}\right)^{p-1}\right](w_2^p - w_1^p)dx.
\end{aligned}
\tag{3.10}
$$

Suppose by contradiction that $D_{R_0}(0) \neq \emptyset$. Thanks to the continuity of $u_i$, $i = 1, 2$, $D_{R_0}(0)$ is an open set. Let $B$ be a small ball such that $\overline{B} \subset D_{R_0}(0)$. As $g$ is nonnegative and $B \subset D_R(\epsilon)$ for all $\epsilon$ small, it follows that

$$
c_p \int_B \frac{|\nabla \ln u_2 - \nabla \ln u_1|^{p+(2-p)^+}}{(|\nabla \ln u_1| + |\nabla \ln u_2|)^{(2-p)^+}}[u_1^p + u_2^p]dx + \int_B b\left[\frac{g(u_2)}{u_2^{p-1}} - \frac{g(u_1)}{u_1^{p-1}}\right](u_2^p - u_1^p)dx \leq J_\epsilon
\tag{3.11}
$$

where $J_\epsilon$ is the third term in (3.10). Observe that, for $\epsilon < 1$, one has in $D_R(\epsilon)$

$$
0 \leq w_2^p - w_1^p \leq p(w_2 - w_1)w_2^{p-1} \leq p(\|u_2\|_{L^\infty(\mathbb{R}^N)} + 1)^p,
$$

which implies

$$
J_\epsilon \leq p\sup_{\mathcal{N}} |a - \beta\Phi|(\|u_2\|_{L^\infty(\mathbb{R}^N)} + 1)^p \int_{D_R(\epsilon)\cap\mathcal{N}} \left|\left(\frac{u_2}{w_2}\right)^{p-1} - \left(\frac{u_1}{w_1}\right)^{p-1}\right| dx.
$$

Notice that $\mathcal{N}$ is bounded and that

$$\lim_{\epsilon \to 0} \left| \left( \frac{u_2}{w_2} \right)^{p-1} - \left( \frac{u_1}{w_1} \right)^{p-1} \right| = 0 \quad \text{pointwise in } \mathbb{R}^N.$$

By Lebesgue dominated convergence theorem, we obtain

(3.12) 
$$\lim_{\epsilon \to 0} J_\epsilon = 0.$$

Hence by letting $\epsilon \to 0$ in (3.11) we obtain

$$\int_B \frac{|\nabla \ln u_2 - \nabla \ln u_1|^{p+(2-p)^+}}{(|\nabla \ln u_1| + |\nabla \ln u_2|)^{(2-p)^+}} [u_1^p + u_2^p] dx = \int_B b \left[ \frac{g(u_2)}{u_2^{p-1}} - \frac{g(u_1)}{u_1^{p-1}} \right] (u_2^p - u_1^p) dx = 0,$$

which, together with condition (A2) on $g$, implies

(3.13) 
$$|\nabla(\ln u_2 - \ln u_1)| \equiv b \equiv 0 \quad \text{in } D_{R_0}(0).$$

Since $b \not\equiv 0$ in $B_{R_0}$, $D_{R_0}(0) \subsetneq B_{R_0}$ and $\partial D_{R_0}(0) \cap B_{R_0} \neq \emptyset$. Hence there is a connected component $\mathcal{O}$ of $D_{R_0}(0)$ such that $\partial D_{R_0}(0) \cap \partial \mathcal{O} \cap B_{R_0} \neq \emptyset$. From (3.13), we deduce that $\ln u_1 - \ln u_2 \equiv$ constant in $\mathcal{O}$, which in turn implies $u_1 = \ell u_2$ in $\mathcal{O}$ for some $\ell > 0$. As $u_1 = u_2$ on $\partial D_{R_0}(0) \cap \partial \mathcal{O} \cap B_{R_0}$, it follows that $\ell = 1$, which contradicts $\mathcal{O} \subset D_{R_0}(0)$. Therefore, we must have $D_{R_0}(0) = \emptyset$ and thus $u_1 \geq u_2$ in $B_{R_0}$. Since $R_0 > 0$ is arbitrarily large, we conclude that $u_1 \geq u_2$ in $\mathbb{R}^N$. □

## 4. Existence/Nonexistence and Uniqueness

### 4.1. Construction of decaying supersolution.

**Proposition 4.1.** Assume $p > 1$ and (A1),(A3) and (A4) are satisfied.

(i) If $0 < \alpha < p$, there exists a bounded, positive, exponentially decaying supersolution of (1.1).
(ii) If $p = \alpha$, there exists a bounded, positive, polynomially decaying supersolution of (1.1).

*Proof.* **Case 1:** $0 < \alpha < p$. Since (A4) holds, $\{x \in \mathbb{R}^N : \Phi(x) < 0\}$ is bounded. Therefore, one can find $R_\Phi > 0$ such that $\Phi \geq 0$ in $B_{R_\Phi}^c$. Set

$$v(x) = Ce^{-\theta|x|^{1-\frac{\alpha}{p}}}, x \in \mathbb{R}^N, \quad \text{with} \quad \theta = \left( \frac{m}{p-1} \right)^{\frac{1}{p}} \frac{p}{p-\alpha}$$

where $C$ will be made precise later. Let $\varepsilon > 0$. Thanks to (A1), there exists $R_\varepsilon > 0$ such that $-a(x) > (m+\varepsilon)|x|^{-\alpha}$ for every $|x| > R_\varepsilon$. Consequently, for $|x| \geq \max\{R_\varepsilon, R_\Phi\}$,

(4.1)
$$-\Delta_p v + (\beta\Phi - a)v^{p-1} \geq |x|^{-\alpha} v^{p-1} \left[ T|x|^{\frac{\alpha}{p}-1} - \theta^p \left( 1 - \frac{\alpha}{p} \right)^p (p-1) + m + \varepsilon \right]$$

where

$$T = \theta^{p-1} \left( 1 - \frac{\alpha}{p} \right)^{p-1} \left( -\alpha + \frac{\alpha}{p} - 1 + N \right).$$

Since $\alpha < p$, we can choose $R_\varepsilon$ large enough so that $|T||x|^{\frac{\alpha}{p}-1} < \varepsilon$ for every $|x| \geq R_\varepsilon$. Therefore, from (4.1), we get, for $|x| \geq \bar{R} := \max\{R_\varepsilon, R_\Phi\}$,

$$-\Delta_p v + (\beta\Phi - a)v^{p-1} \geq |x|^{-\alpha} v^{p-1} \left[ m - \theta^p \left( 1 - \frac{\alpha}{p} \right)^p (p-1) \right] = 0.$$

We next make use of hypothesis (A3). Define $\overline{v} = s_0 \chi_{\overline{B}_R} + v\chi_{B_R^c}$ where $s_0$ is the positive constant given in (A3). Choose $C = s_0 e^{\theta \bar{R}^{1-\frac{\alpha}{p}}}$ then $\overline{v}$ is a positive, decaying weak supersolution of (1.1) in $\mathbb{R}^N$.

**Case 2:** $\alpha = p$. Let $\varepsilon > 0$. By a similar argument as in Case 1, we can construct a positive, decaying super solution of (1.1) under the form

$$\overline{v}(x) = s_0 \chi_{\overline{B}_R}(x) + C|x|^{-\tilde{\theta}} \chi_{B_R^c}(x), \quad x \in \mathbb{R}^N$$

where $\bar{R} = \bar{R}(\varepsilon, \Phi)$, $\tilde{\theta}$ is the unique positive solution of the algebraic equation

$$(p-1)\gamma^p - (N-p)\gamma^{p-1} - m = 0$$

and $C = s_0 \bar{R}^{\tilde{\theta}}$. $\qquad\qquad\qquad\qquad\qquad\qquad\qquad\qquad\qquad\qquad \square$

4.2. **A-priori estimates in the case $p \geq 2$.** In this subsection, we establish a-priori estimates for positive solutions of (1.1) which are crucial to obtain the nonexistence and uniqueness results. Decay phenomena are different between two cases $\alpha \in [0, p)$ and $\alpha = p$.

**Proposition 4.2.** Let $p \geq 2$, $\alpha \in [0, p)$, $m > 0$ and $u$ be a positive function satisfying

$$(4.2) \qquad \liminf_{|x| \to \infty} |x|^\alpha \left( \frac{\Delta_p u}{u^{p-1}} - \frac{m}{|x|^\alpha} \right) > 0 \qquad \text{and} \qquad \limsup_{|x| \to \infty} \frac{u(x)}{e^{\overline{\omega}|x|^{1-\frac{\alpha}{p}}}} < \infty,$$

where $\overline{\omega}$ is given by (1.8). Then

$$(4.3) \qquad\qquad\qquad\qquad \lim_{|x| \to \infty} e^{\overline{\omega}|x|^{1-\frac{\alpha}{p}}} u(x) = 0.$$

*Proof.* By (4.2), for any $\varepsilon > 0$ there exists $R = R(\varepsilon)$ such that

$$\Delta_p u \geq (m+\varepsilon)|x|^{-\alpha} u^{p-1} \qquad \text{for } |x| \geq R.$$

Set $\mathcal{L}_\varepsilon[\phi] := -\Delta_p \phi + (m+\varepsilon)|x|^{-\alpha}\phi^{p-1}$. It is easy to see that $\mathcal{L}_\varepsilon[u] \leq 0$ in $B_{\bar{R}}^c$. For any $\rho > 0$, set

$$w_\rho^1(x) = e^{(R+\rho)^{1-\frac{\alpha}{p}}(\tau-\omega)} e^{\omega|x|^{1-\frac{\alpha}{p}}}, \qquad w_\rho^2(x) = e^{R^{1-\frac{\alpha}{p}}(\tau+\omega)} e^{-\omega|x|^{1-\frac{\alpha}{p}}},$$

$$w_\rho = w_\rho^1 + w_\rho^2$$

where $\omega, \tau$ and $R$ will be chosen later.

We will estimate $\Delta_p w_\rho$. Observe that, in $\mathbb{R}^N \setminus \{0\}$,

$$(4.4) \qquad \Delta_p w_\rho = (p-2)|\nabla w_\rho|^{p-4} \left\langle D^2 w_\rho \nabla w_\rho, \nabla w_\rho \right\rangle + |\nabla w_\rho|^{p-2} \Delta w_\rho.$$

By Cauchy-Schwarz inequality, one has

$$(4.5) \qquad \Delta_p w_\rho \leq \left( (p-2)N \max_{ij} |\partial_{ij} w_\rho| + |\Delta w_\rho| \right) |\nabla w_\rho|^{p-2}.$$

Next, we look for an upper bound for the right-hand side of (4.5). Direct computation yields, for every $x \neq 0$,

$$\nabla w_\rho = \omega \left( 1 - \frac{\alpha}{p} \right) |x|^{-\frac{\alpha}{p} - 1} x \, w_\rho^1 - \omega \left( 1 - \frac{\alpha}{p} \right) |x|^{-\frac{\alpha}{p} - 1} x \, w_\rho^2,$$

thus

$$(4.6) \qquad |\nabla w_\rho|^{p-2} \leq \omega^{p-2} \left( 1 - \frac{\alpha}{p} \right)^{p-2} |x|^{-\frac{\alpha(p-2)}{p}} w_\rho^{p-2}.$$

For every $1 \leq i, j \leq N$ and $x \neq 0$,

$$\begin{aligned}
\partial_{ij} w_\rho &= \omega^2 \left( 1 - \frac{\alpha}{p} \right)^2 \frac{x_i x_j}{|x|^{2 + \frac{2\alpha}{p}}} w_\rho^1 - \omega \left( 1 - \frac{\alpha^2}{p^2} \right) \frac{x_i x_j}{|x|^{3 + \frac{\alpha}{p}}} w_\rho^1 \\
&+ \delta(i-j)\omega \left( 1 - \frac{\alpha}{p} \right) \frac{w_\rho^1}{|x|^{1 + \frac{\alpha}{p}}} + \omega^2 \left( 1 - \frac{\alpha}{p} \right)^2 \frac{x_i x_j}{|x|^{2 + \frac{2\alpha}{p}}} w_\rho^2 \\
&+ \omega \left( 1 - \frac{\alpha^2}{p^2} \right) \frac{x_i x_j}{|x|^{3 + \frac{\alpha}{p}}} w_\rho^2 - \delta(i-j)\omega \left( 1 - \frac{\alpha}{p} \right) \frac{w_\rho^2}{|x|^{1 + \frac{\alpha}{p}}}
\end{aligned}$$

where $\delta$ is the Dirac function. Since $|x_i x_j| \leq |x|^2$, it follows that

$$(4.7) \qquad \begin{aligned}
|\partial_{ij} w_\rho| &\leq \omega^2 \left( 1 - \frac{\alpha}{p} \right)^2 \frac{w_\rho}{|x|^{\frac{2\alpha}{p}}} + \omega \left( 1 - \frac{\alpha^2}{p^2} \right) \frac{w_\rho}{|x|^{1 + \frac{\alpha}{p}}} + \omega \left( 1 - \frac{\alpha}{p} \right) \frac{w_\rho}{|x|^{1 + \frac{\alpha}{p}}} \\
&\leq \left[ \omega^2 \left( 1 - \frac{\alpha}{p} \right)^2 |x|^{1 - \frac{\alpha}{p}} + \omega \left( 1 - \frac{\alpha^2}{p^2} \right) + \omega \left( 1 - \frac{\alpha}{p} \right) \right] |x|^{-1 - \frac{\alpha}{p}} w_\rho.
\end{aligned}$$

Consequently,

$$(4.8) \quad |\Delta w_\rho| \leq N \left[ \omega^2 \left( 1 - \frac{\alpha}{p} \right)^2 |x|^{1 - \frac{\alpha}{p}} + \omega \left( 1 - \frac{\alpha^2}{p^2} \right) + \omega \left( 1 - \frac{\alpha}{p} \right) \right] |x|^{-1 - \frac{\alpha}{p}} w_\rho.$$

Combining (4.5)-(4.8), we have

$$\Delta_p w_\rho \leq N(p-1)\omega^{p-1} \left( 1 - \frac{\alpha}{p} \right)^{p-1} \left[ \omega \left( 1 - \frac{\alpha}{p} \right) |x|^{1 - \frac{\alpha}{p}} + \frac{2p + \alpha}{p} \right] |x|^{-\alpha - 1 + \frac{\alpha}{p}} w_\rho^{p-1}.$$

Put

$$A = N(p-1)\omega^{p-1}\left(1-\frac{\alpha}{p}\right)^{p-1}\frac{2p+\alpha}{p}.$$

As $|x| \geq R$, one gets

$$(4.9) \qquad \mathcal{L}_\varepsilon[w_\rho] \geq |x|^{-\alpha}w_\rho^{p-1}\left[-N(p-1)\omega^p\left(1-\frac{\alpha}{p}\right)^p - A|x|^{\frac{\alpha}{p}-1} + m + \varepsilon\right].$$

One can choose $R$ and $\omega$ such that the right-hand side of (4.9) is nonnegative. Indeed, since $\alpha < p$, $|x|^{\frac{\alpha}{p}-1} \to 0$ as $|x| \to \infty$, hence there exists $R(\varepsilon)$ such that, for every $R > R(\varepsilon)$,

$$A|x|^{\frac{\alpha}{p}-1} \leq \frac{\varepsilon}{2} \quad \forall x \in B_R^c.$$

Put

$$\omega := \left(\frac{2m+\varepsilon}{2N(p-1)}\right)^{1/p}\frac{p}{p-\alpha}.$$

With such $R$ and $\omega$, we obtain $\mathcal{L}_\varepsilon[w_\rho] \geq 0$ in $B_{R+\rho} \setminus B_\rho$.

We next show that $w_\rho$ dominates $u$ on $\partial B_{R+\rho} \cup \partial B_\rho$. Indeed, by (4.2), one can finds $C > 0$ such that $u(x) \leq Ce^{\overline{\omega}|x|^{1-\frac{\alpha}{p}}}$ in $\mathbb{R}^N$. Therefore, we can take $\tau$ arbitrarily in $(\overline{\omega}, \omega)$ and $R$ sufficiently large such that for any $\rho > 0$, one has

$$\begin{cases} w_\rho(x) \geq e^{R^{1-\frac{\alpha}{p}}\tau} \geq Ce^{R^{1-\frac{\alpha}{p}}\overline{\omega}} \geq u(x), & \text{as } |x| = R \\ w_\rho(x) \geq e^{(R+\rho)^{1-\frac{\alpha}{p}}\tau} \geq Ce^{(R+\rho)^{1-\frac{\alpha}{p}}\overline{\omega}} \geq u(x), & \text{as } |x| = R+\rho. \end{cases}$$

Fix such $\omega, \tau$ and $R$. Applying the weak comparison principle [21], we obtain

$$u(x) \leq w_\rho(x) = e^{(R+\rho)^{1-\frac{\alpha}{p}}(\tau-\omega)}e^{\omega|x|^{1-\frac{\alpha}{p}}} + e^{R^{1-\frac{\alpha}{p}}(\tau+\omega)}e^{-\omega|x|^{1-\frac{\alpha}{p}}} \qquad \text{in } B_{R+\rho} \setminus B_R.$$

Sending $\rho \to \infty$ yields

$$u(x) \leq e^{R^{1-\frac{\alpha}{p}}(\tau+\omega)}e^{-\omega|x|^{1-\frac{\alpha}{p}}} \qquad \text{in } \mathbb{R}^N \setminus B_R.$$

The fact $\omega > \overline{\omega}$ confirms the proof. $\qquad \square$

When $\alpha = p$ we obtain the following result.

**Proposition 4.3.** Let $p \geq 2$, $\alpha = p$, $m > 0$ and $u$ be a positive function satisfying

$$(4.10) \qquad \liminf_{|x|\to\infty}|x|^p\left(\frac{\Delta_p u}{u^{p-1}} - \frac{m}{|x|^p}\right) > 0 \qquad \text{and} \qquad \limsup_{|x|\to\infty}\frac{u(x)}{|x|^{\gamma_0}} < \infty,$$

where $\gamma_0$ is the unique positive solution of (1.9). Then

$$(4.11) \qquad \lim_{|x|\to\infty}|x|^{\gamma_0}u(x) = 0.$$

*Proof.* By (4.10), for any $\varepsilon > 0$, there exists $R = R(\varepsilon)$ such that

$$\Delta_p u \geq (m + \varepsilon)|x|^{-p} u^{p-1} \qquad \text{for } |x| \geq R.$$

Set $\mathcal{L}_\varepsilon[\phi] = -\Delta_p \phi + (m + \varepsilon)|x|^{-p} \phi^{p-1}$. Obviously $\mathcal{L}_\varepsilon[u] \leq 0$ in $\mathbb{R}^N \setminus B_R$. For every $\rho > 0$, we look for a supersolution of $\mathcal{L}_\varepsilon[\phi] = 0$ in $B_{R+\rho} \setminus B_R$. For $\rho > 0$, set

$$w_\rho(x) := \frac{C_1}{\ln(R + \rho)}|x|^\gamma + C_2 |x|^{-\gamma}, \quad x \in \mathbb{R}^N \setminus \{0\}$$

where $C_1, C_2, \gamma > 0$. One can proceed as in the proof of Proposition 4.2 with some obvious modifications to show that $w_\rho$ is a supersolution of $\mathcal{L}_\varepsilon[\phi] = 0$ in $\mathbb{R}^N \setminus B_R$. However, for the sake of completion, we present the detailed computations. It is easy to get, for $x \neq 0$,

$$\nabla w_\rho = \frac{C_1}{\ln(R + \rho)}\gamma |x|^{\gamma-2} x - C_2 \gamma |x|^{-\gamma-2} x,$$

which implies

$$(4.12) \qquad |\nabla w_\rho|^{p-2} \leq \gamma^{p-2}|x|^{-(p-2)} w_\rho^{p-2}.$$

For $1 \leq i, j \leq N$ and $x \neq 0$,

$$(4.13) \qquad \begin{aligned} \partial_{ij} w_\rho &= \frac{C_1}{\ln(R + \rho)} \left[ \gamma(\gamma - 2)|x|^{\gamma-4} x_i x_j + \delta(i - j)\gamma |x|^{\gamma-2} \right] \\ &\quad + C_2 \left[ \gamma(\gamma + 2)|x|^{-\gamma-4} x_i x_j - \delta(i - j)\gamma |x|^{-\gamma-2} \right]. \end{aligned}$$

Hence

$$\begin{aligned} \Delta w_\rho &= \frac{C_1}{\ln(R + \rho)} \left[ \gamma(\gamma - 1)|x|^{\gamma-2} + \gamma(N - 1)|x|^{\gamma-2} \right] \\ &\quad + C_2 \left[ \gamma(\gamma + 1)|x|^{-\gamma-2} - \gamma(N - 1)|x|^{-\gamma-2} \right], \end{aligned}$$

and

$$(4.14) \qquad |\Delta w_\rho| \leq \gamma(\gamma + N)|x|^{-2} w_\rho.$$

On the other hand, by the inequality $|x_i x_j| \leq |x|^2$, we deduce from (4.13) that, for $1 \leq i, j \leq N$ and $x \neq 0$,

$$(4.15) \qquad |\partial_{ij} w_\rho| \leq \gamma(\gamma + 3)|x|^{-2} w_\rho.$$

Combining (4.5) and (4.12)-(4.15), we get

$$\Delta_p w_\rho \leq \left[ ((p - 2)N + 1)\gamma^p + (3p - 5)N\gamma^{p-1} \right] |x|^{-p} w_\rho^{p-1}.$$

Therefore, as $|x| \geq R$, we have

$$(4.16) \qquad \mathcal{L}_\varepsilon[w_\rho] \geq \left[ -((p - 2)N + 1)\gamma^p - (3p - 5)N\gamma^{p-1} + m + \varepsilon \right] |x|^{-p} w_\rho^{p-1}.$$

Since $p \geq 2$, there exists a unique positive solution $\gamma_\varepsilon$ of the equation

$$((p - 2)N + 1)\gamma^p + (3p - 5)N\gamma^{p-1} - m - \varepsilon = 0.$$

Obviously $\gamma_\varepsilon > \gamma_0$ where $\gamma_0$ is the unique positive solution of (1.9). By taking $\gamma = \gamma_\varepsilon$ in (4.16), we get $\mathcal{L}_\varepsilon[w_\rho] \geq 0$ in $B_{R+\rho} \setminus B_R$ for all $\rho > 0$.

We next show that $C_1$ and $C_2$ can be chosen large enough, independent of $\rho$, such that $u \leq w_\rho$ on $\partial B_{R+\rho} \cup B_R$. Indeed, it follows from (1.11) that there exists $C > 0$ such that $u(x) \leq C|x|^{\gamma_0}$ for $|x| \geq R$. Choosing $C_1 = C$, for $\rho$ large, one has

$$u(x) \leq C|x|^{\gamma_0} \leq C\frac{|x|^{\gamma_\varepsilon - \gamma_0}}{\ln(R+\rho)}|x|^{\gamma_0} \leq w_\rho(x) \qquad x \in \partial B_{R+\rho}.$$

We choose $C_2$ large enough, depending on $R$, such that $u \leq w_\rho$ on $\partial B_R$.

By the weak comparison principle, we obtain

$$u(x) \leq w_\rho(x) = \frac{C_1}{\ln(R+\rho)}|x|^{\gamma_\varepsilon} + C_2|x|^{-\gamma_\varepsilon} \qquad x \in B_{R+\rho} \setminus B_R.$$

Since $C_1$ are $C_2$ are independent of $\rho$, by letting $\rho \to \infty$, we finally derive

$$u(x) \leq C_2|x|^{-\gamma_\varepsilon} \qquad x \in \mathbb{R}^N \setminus B_R.$$

The proof is complete since $\gamma_\varepsilon > \gamma_0$. $\qquad\qquad\qquad\qquad\qquad\qquad\qquad\square$

4.3. **Proof of Theorem 1.3 and Theorem 1.4.** This section is devoted to demonstration of the main results of the paper.

**Proof of Theorem 1.3.**

**Case 1:** $\lambda_\beta < 0$. We use sub and super solutions argument to prove the existence of positive solutions of (1.1).

We first construct a subsolution. By Theorem 1.2, there exists $R_\beta > 0$ such that $\lambda(\mathcal{K}_\beta, B_{R_\beta}) < 0$. Let $\varphi$ be the positive eigenfunction associated with $\lambda(\mathcal{K}_\beta, B_{R_\beta})$, normalized by $\varphi(0) = 1$. Set $\underline{u}_\beta(x) = \delta\varphi$ where $\delta > 0$ is chosen later on and $M_\beta := \max_{B_{R_\beta}} \varphi$. In $B_{R_\beta}$, we get

$$
\begin{aligned}
-\Delta_p \underline{u}_\beta + (\beta\Phi - a)\underline{u}_\beta^{p-1} + bg(\underline{u}_\beta) \quad &\leq (\delta\varphi)^{p-1}\left(\lambda(\mathcal{K}_\beta, B_{R_\beta}) + \|b\|_{L^\infty(\mathbb{R}^N)} \frac{g(\delta\varphi)}{(\delta\varphi)^{p-1}}\right) \\
&\leq (\delta\varphi)^{p-1}\left(\lambda(\mathcal{K}_\beta, B_{R_\beta}) + \|b\|_{L^\infty(\mathbb{R}^N)} \frac{g(\delta M_\beta)}{(\delta M_\beta)^{p-1}}\right).
\end{aligned}
\tag{4.17}
$$

By (A2), one can choose $\delta$ small, depending on $\beta$, in such a way that $\delta M_\beta < s_0$ (where $s_0$ is the positive constant in (A3)) and

$$\frac{g(\delta M_\beta)}{(\delta M_\beta)^{p-1}} < -\frac{\lambda(\mathcal{K}_\beta, B_{R_\beta})}{2\|b\|_{L^\infty(\mathbb{R}^N)}}.\tag{4.18}$$

Consequently the right-hand side of (4.17) is negative. Put

$$\underline{U}_\beta(x) = \begin{cases} \chi_{B_{R_\beta}}(x)\underline{u}_\beta(x) & \text{if } x \in B_{R_\beta} \\ 0 & \text{otherwise} \end{cases}$$

then $\underline{U}_\beta$ is a nonnegative weak subsolution of (1.1) in $\mathbb{R}^N$ and $\underline{U}_\beta < s_0$.

For any $R > R_\beta$, $\underline{U}_\beta$ and $s_0$ are respectively sub and super solutions of

$$(4.19) \qquad \begin{cases} -\Delta_p u + (\beta\Phi - a)u + bg(u) & = 0 \quad \text{in } B_R \\ u & = 0 \quad \text{on } \partial B_R. \end{cases}$$

By applying the sub and super solution theorem (see for instance [24, Theorem 3.1]), we derive that there exists a weak solution $u_{\beta,R}$ of (4.19) in $B_R$ such that $\underline{U}_\beta \le u_{\beta,R} \le s_0$ in $B_R$. By local regularity for quasilinear elliptic equations (see [19]) and standard argument, $\{u_{\beta,R}\}$ converges, as $R \to \infty$, in $C^1_{\text{loc}}(\mathbb{R}^N)$ to a function $u_\beta$ which is a weak solution of (1.1) in $\mathbb{R}^N$ and satisfies $\underline{U}_\beta \le u_\beta \le s_0$. Since $u_\beta(0) \ge \underline{U}_\beta(0) = \delta > 0$, by Harnack inequality (see [29], [31]), we obtain $u_\beta > 0$ in $\mathbb{R}^N$. Since $u_\beta \le s_0$, it follows that $u_\beta \in \mathcal{S}_p$. From Proposition 4.2 and Proposition 4.3, we derive (1.10) if $\alpha \in [0,p)$ and (1.11) if $\alpha = p$.

The uniqueness is a direct consequence of decay property and comparison principle Theorem 3.1.

**Case 2: $\lambda_\beta \ge 0$.** Suppose by contradiction that there exists a weak solution $u$ of (1.1) belonging to $\mathcal{S}_p$. Thanks to (A1), Proposition 4.2 and Proposition 4.3 imply that $u$ satisfies either (1.10) or (1.11); in particular, $u$ decays. Let $\varphi$ be a positive weak solution of

$$(4.20) \qquad \mathcal{K}_\beta[\varphi] = \lambda_\beta \varphi^{p-1}$$

in $\mathbb{R}^N$, normalized by $\varphi(0) < u(0)$ (the existence of $\varphi$ is guaranteed by Theorem 1.2 (ii)). Since $\lambda_\beta \ge 0$, it follows that $\varphi$ is a positive supersolution of (1.1). By Theorem 3.1, we have $u \le \varphi$ in $\mathbb{R}^N$. This contradiction completes the proof. $\qquad \square$

**Proof of Theorem 1.4.**
**Case 1: $\lambda_\beta < 0$.** Let $\overline{v}$ be the decaying positive supersolution of (1.1) in $\mathbb{R}^N$ constructed in the proof of Proposition 4.1. Let $\underline{U}_\beta$ be the subsolution of (1.1) constructed in Case 1 of the proof of Theorem 1.3 with $\delta$ being chosen such that $\delta M_\beta < \inf_{B_{R_\beta}} \overline{v}$ and (4.18) holds. Then $\overline{U}_\beta < v$ in $\mathbb{R}^N$. The rest of the proof can be proceeded as in Case 1 of the proof of Theorem 1.3 and we omit it.

**Case 2: $\lambda_\beta \ge 0$.** Suppose by contradiction that there is a solution $u$ of (1.1) belonging to $\mathcal{S}_p$. Thanks to Theorem 3.1, by a similar argument as in Case 2 of the proof of Theorem 1.3, we derive a contradiction. $\qquad \square$

## 5. THRESHOLD VALUE AND ASYMPTOTIC BEHAVIORS

Throughout this section, we assume that $\Phi \ge 0$, $\Phi \in L^\infty_{\text{loc}}(\mathbb{R}^N)$. For $\beta \ge 0$, let $\lambda_\beta$ be the generalized principal eigenvalue of $\mathcal{K}_\beta$ in $\mathbb{R}^N$, i.e. $\lambda_\beta = \lambda(\mathcal{K}_\beta, \mathbb{R}^N)$. Applying Theorem 1.2, (iii) with $\Omega = \mathbb{R}^N$ and $V = \beta\Phi - a$, one gets

$$(5.1) \qquad \lambda_\beta = \inf_{\phi \in C^1_c(\mathbb{R}^N)\backslash\{0\}} \frac{\int_{\mathbb{R}^N}(|\nabla\phi|^p + (\beta\Phi - a)|\phi|^p)dx}{\int_{\mathbb{R}^N} |\phi|^p dx}.$$

We recall an important hypothesis in this section:

(A5)     $\Phi(0) \geq 0$, $\Phi > 0$ in $\mathbb{R}^N \setminus \{0\}$, $\Phi \in L_{\text{loc}}^\infty(\mathbb{R}^N)$, $\lim_{|x|\to\infty} \Phi(x) = \infty$ and
         the embedding $W_\Phi^{1,p}(\mathbb{R}^N) \hookrightarrow L^p(\mathbb{R}^N)$ is compact.

The next result shows that the infimum can be taken over $W_\Phi^{1,p}(\mathbb{R}^N) \setminus \{0\}$ and the the corresponding variational problem admits a minimizer.

**Proposition 5.1.** Assume that $p > 1$, $\beta > 0$ and (A5) is satisfied. There holds

$$(5.2) \qquad \lambda_\beta = \min_{\phi \in W_\Phi^{1,p}(\mathbb{R}^N) \setminus \{0\}} \frac{\int_{\mathbb{R}^N}(|\nabla \phi|^p + (\beta\Phi - a)|\phi|^p)dx}{\int_{\mathbb{R}^N}|\phi|^p dx}.$$

and $\lambda_\beta$ is achieved at some positive function $\varphi_\beta \in W_\Phi^{1,p}(\mathbb{R}^N)$.

*Proof.* Put

$$\hat{\lambda}_\beta := \inf_{\phi \in W_\Phi^{1,p}(\mathbb{R}^N) \setminus \{0\}} \frac{\int_{\mathbb{R}^N}(|\nabla \phi|^p + (\beta\Phi - a)|\phi|^p)dx}{\int_{\mathbb{R}^N}|\phi|^p dx}.$$

We will show that $\lambda_\beta = \hat{\lambda}_\beta$. Since $C_c^1(\mathbb{R}^N) \subset W_\Phi^{1,p}(\mathbb{R}^N)$, from (5.1), one gets $\lambda_\beta \geq \hat{\lambda}_\beta \geq -\sup_{\mathbb{R}^N} a$. We next prove the inverse inequality. By the definition of $\hat{\lambda}_\beta$, for every $\delta > 0$, there exists $\phi_\delta \in W_\Phi^{1,p}(\mathbb{R}^N)$, $\phi_\delta \neq 0$ such that

$$(5.3) \qquad \hat{\lambda}_\beta + \delta \geq \frac{\int_{\mathbb{R}^N}(|\nabla \phi_\delta|^p + (\beta\Phi - a)|\phi_\delta|^p)dx}{\int_{\mathbb{R}^N}|\phi_\delta|^p dx}.$$

Let $\{\eta_n\}_{n>1} \subset C^1(\mathbb{R}^N)$ be a sequence of functions such that $0 \leq \eta_n \leq 1$, $\eta_n = 0$ in $B_n^c$, $\eta_n = 1$ in $B_{n-1}$ and $\|\nabla \eta_n\|_{L^\infty(\mathbb{R}^N)} \leq M$ for every $n > 1$. Set $\phi_{\delta,n} = \eta_n \phi_\delta$. Since $\phi_{\delta,n} \in W_0^{1,p}(B_n)$, thanks to (1.4), we get

$$(5.4) \qquad \frac{\int_{\mathbb{R}^N}(|\nabla \phi_{\delta,n}|^p + (\beta\Phi - a)|\phi_{\delta,n}|^p)dx}{\int_{\mathbb{R}^N}|\phi_{\delta,n}|^p dx} \geq \lambda(\mathcal{K}_\beta, B_n).$$

On the other hand,

$$|\nabla \phi_{\delta,n}| = |\eta_n \nabla \phi_\delta + \phi_\delta \nabla \eta_n| \leq |\nabla \phi_\delta| + M|\phi_\delta|,$$

from which it follows

$$(5.5) \qquad |\nabla \phi_{\delta,n}|^p + \beta\Phi|\phi_{\delta,n}|^p \leq 2^{p-1}(|\nabla \phi_\delta|^p + \beta\Phi|\phi_\delta|^p + M^p|\phi_\delta|^p).$$

Since $W_\Phi^{1,p}(\mathbb{R}^N) \hookrightarrow L^p(\mathbb{R}^N)$, $\phi_\delta \in L^p(\mathbb{R}^N)$, whence the right hand-side of (5.5) belongs to $L^1(\mathbb{R}^N)$. The Lebesgue dominated convergence theorem, along with (5.3) and (5.4), implies that

$$\hat{\lambda}_\beta + \delta \geq \frac{\int_{\mathbb{R}^N}(|\nabla \phi_\delta|^p + (\beta\Phi - a)|\phi_\delta|^p)dx}{\int_{\mathbb{R}^N}|\phi_\delta|^p dx} \geq \lim_{R\to\infty}\lambda(\mathcal{K}_\beta, B_n) = \lambda_\beta.$$

Letting $\delta \to 0$ yields $\hat{\lambda}_\beta \geq \lambda_\beta$, thus $\hat{\lambda}_\beta = \lambda_\beta$.

Due to the compact embedding $W_\Phi^{1,p}(\mathbb{R}^N) \hookrightarrow L^p(\mathbb{R}^N)$, it is classical that $\lambda_\beta$ is achieved at a function $\varphi_\beta \in W_\Phi^{1,p}(\mathbb{R}^N)$ with $\|\varphi_\beta\|_{L^p(\mathbb{R}^N)} = 1$. The positivity of $\varphi_\beta$ can be obtained by using a similar argument as in [22, Proposition 5.3]. $\quad\square$

**Proposition 5.2.** Assume $p > 1$ and (A5) is satisfied. Then the mapping $\beta \mapsto \lambda_\beta$ is continuous, increasing and concave on $(0, \infty)$. Moreover, $\lim_{\beta \to 0} \lambda_\beta = \lambda_0$ and $\lambda_\beta > 0$ for $\beta$ large enough.

*Proof. Step 1: The mapping $\beta \mapsto \lambda_\beta$ is continuous, increasing and concave on $(0, \infty)$.* Take arbitrarily $\beta > 0$ and $\beta' > 0$. From (5.2), we get

$$\lambda_{\beta+\beta'} \leq \int_{\mathbb{R}^N} [|\nabla\varphi_\beta|^p + ((\beta+\beta')\Phi - a)\varphi_\beta^p]dx = \lambda_\beta + \beta'\int_{\mathbb{R}^N}\Phi\varphi_\beta^p dx.$$

Similarly,

$$\lambda_\beta \leq \lambda_{\beta+\beta'} - \beta'\int_{\mathbb{R}^N}\Phi\varphi_{\beta+\beta'}^p dx.$$

Therefore

$$0 < \beta'\int_{\mathbb{R}^N}\Phi\varphi_{\beta+\beta'}^p dx \leq \lambda_{\beta+\beta'} - \lambda_\beta \leq \beta'\int_{\mathbb{R}^N}\Phi\varphi_\beta^p dx.$$

Consequently, $\beta \mapsto \lambda_\beta$ is increasing and locally Lipschitz on $(0, \infty)$.

For every $t \in [0, 1]$ and $\beta > 0$, $\beta' > 0$, by (5.2), it is easy to see that

$$\lambda_{t\beta+(1-t)\beta'} \geq t\lambda_\beta + (1-t)\lambda_{\beta'}.$$

Thus $\beta \mapsto \lambda_\beta$ is concave on $(0, \infty)$.

*Step 2: We prove that $\lim_{\beta \to 0} \lambda_\beta = \lambda_0$.* Notice that $\lambda_{\beta'} \geq \lambda_\beta \geq \lambda_0$ for every $0 < \beta < \beta'$ since $\Phi \geq 0$. Therefore $\underline{\lambda} := \lim_{\beta \to 0} \lambda_\beta \geq \lambda_0$. By Theorem 1.2, for any $\delta > 0$, there exists $\phi_\delta \in C_c^1(\mathbb{R}^N)$ with $\|\phi_\delta\|_{L^p(\mathbb{R}^N)} = 1$ such that

$$\lambda_0 + \delta \geq \int_{\mathbb{R}^N} (|\nabla\phi_\delta|^p - a|\phi_\delta|^p)dx.$$

Since $\Phi \in L_{\text{loc}}^\infty(\mathbb{R}^N)$, $\phi_\delta \in C_c^1(\mathbb{R}^N)$, we have $\Phi|\phi_\delta|^p \in L^1(\mathbb{R}^N)$. Choosing $\beta$ such that $\beta\int_{\mathbb{R}^N}\Phi|\phi_\delta|^p dx < \delta$, we get

$$\lambda_0 + 2\delta \geq \int_{\mathbb{R}^N} [|\nabla\phi_\delta|^p + (\beta\Phi - a)|\phi_\delta|^p]dx \geq \lambda_\beta \geq \underline{\lambda}.$$

Letting $\delta \to 0$ yields $\lambda_0 \geq \underline{\lambda}$, which leads to $\lambda_0 = \underline{\lambda}$.

*Step 3: We show that $\lambda_\beta > 0$ for $\beta$ large.* Suppose by contradiction that for every $\beta \in (0, \infty)$, $\lambda_\beta \leq 0$. It follows that

$$0 \geq \int_{\mathbb{R}^N} (|\nabla\varphi_\beta|^p + (\beta\Phi - a)\varphi_\beta^p)dx \geq \beta\int_{\mathbb{R}^N}\Phi\varphi_\beta^p dx - \sup_{\mathbb{R}^N} a.$$

This in turn implies

$$\int_{\mathbb{R}^N} \Phi \varphi_\beta^p dx \leq \frac{1}{\beta} \sup_{\mathbb{R}^N} a.$$

Therefore $\varphi_\beta \to 0$ in $L^p(\mathbb{R}^N \setminus B_R)$ as $\beta \to \infty$ for every $R > 0$. On the other hand, since $\{\varphi_\beta\}$ is bounded in $W^{1,p}(\mathbb{R}^N)$ and the embedding $W^{1,p}_\Phi(\mathbb{R}^N) \hookrightarrow L^p(\mathbb{R}^N)$ is compact (hypothesis (A5)), we deduce that up to a subsequence, $\{\varphi_\beta\}$ converges strongly in $L^p(\mathbb{R}^N)$. Therefore $\varphi_\beta \to 0$ in $L^p(\mathbb{R}^N)$. This is a contradiction since $\|\varphi_\beta\|_{L^p(\mathbb{R}^N)} = 1$ for every $\beta > 0$.                    □

**Theorem 5.3.** Assume $p > 1$ and (A5) holds. If $\lambda_0 < 0$ then there exists $\beta^*$ such that $\lambda_{\beta^*} = 0$, $\lambda_\beta < 0$ for every $\beta < \beta^*$ and $\lambda_\beta > 0$ for every $\beta > \beta^*$.

*Proof.* Theorem 5.3 is a direct consequence of Proposition 5.2.                    □

Under assumptions (A1), (A2), (A3) and (A5), if $\lambda_0 < 0$, by Theorem 5.3, Theorem 1.3 and Theorem 1.4 we deduce that for each $\beta \in [0, \beta^*)$ there exists a unique solution $u_\beta$ of (1.1) in $\mathcal{S}_p$. Some qualitative properties of the sequence $\{u_\beta\}$ are presented in the following result.

**Proposition 5.4.** Assume $p > 1$, (A1), (A2), (A3) and (A5) hold. If $\lambda_0 < 0$ then for all $0 < q \leq \infty$ there hold

$$(5.6) \qquad \lim_{\beta \to 0} \|u_\beta - u_0\|_{L^q(\mathbb{R}^N)} = 0 \qquad \text{and} \qquad \lim_{\beta \to \beta^*} \|u_\beta\|_{L^q(\mathbb{R}^N)} = 0.$$

*Proof.* Since $\Phi \geq 0$, by Theorem 3.1, we deduce that $\{u_\beta\}$ is nonincreasing with respect to $\beta$. By local regularity for quasilinear elliptic equations [19, 23], $\{u_\beta\}$ converges, as $\beta \to 0$, in $C^1_{\text{loc}}(\mathbb{R}^N)$ to a function $\hat{u}_0$ which is a weak solution of

$$(5.7) \qquad -\Delta_p \hat{u}_0 - a\hat{u}_0^{p-1} + bg(\hat{u}_0) = 0 \quad \text{in } \mathbb{R}^N.$$

Since $u_\beta \leq u_0$ for every $\beta > 0$, it follows that $\hat{u}_0 \leq u_0$. Therefore $\hat{u}_0$ is a decaying positive solution of (5.7). By Theorem 1.3 and Theorem 1.4, $\hat{u}_0 \equiv u_0$.

Similarly, $u_\beta$ converges as $\beta \to \beta^*$ in $C^1_{\text{loc}}(\mathbb{R}^N)$ to $u_{\beta^*}$ which is a weak solution of

$$(5.8) \qquad -\Delta_p u_{\beta^*} + (\beta^* \Phi - a)u_{\beta^*}^{p-1} + bg(u_{\beta^*}) = 0 \quad \text{in } \mathbb{R}^N.$$

Since $\lambda_{\beta^*} = 0$, by Theorem 1.3, Theorem 1.4, (5.8) admits no decaying positive solution. As $0 \leq u_{\beta^*} \leq u_0$, it follows that $u_{\beta^*} \equiv 0$.

Finally (5.6) follows from the monotone convergence theorem.                    □

**Proof of Theorem 1.5.** By combining Theorem 1.3, Theorem 1.4, Theorem 5.3 and Proposition 5.4, we obtain desired results easily.                    □

**Acknowledgements.** The first author was supported by the Israel Science Foundation founded by the Israel Academy of Sciences and Humanities, through grant 91/10 and partially supported by the Technion Postdoctoral fellowship. The second author is supported by PhD fellowship within the framework of ERC Grant (FP7/2007-2013), agreement number 321186: "Reaction-Diffusion Equations, Propagation and Modelling" held by Henri Berestycki. We would like to thank Professor Y. Pinchover for interesting discussions. Lastly, we thank the anonymous referees for their helpful comments.

Department of Mathematics, Technion, Haifa 32000, ISRAEL and Ho Chi Minh City University of Pedagogy, Ho Chi Minh City, Vietnam
*E-mail address*: nguyenphuoctai.hcmup@gmail.com

CAMS - École des Hautes Études en Sciences Sociales, 190-198 avenue de France, 75013, Paris, France
*E-mail address*: vhhungkhtn@gmail.com